# Sensitivity Factors based Transmission Network Topology Control for Violation Relief

Xingpeng Li[1*], Akshay Korad[2], Pranavamoorthy Balasubramanian[3]

[1] Department of Electrical and Computer Engineering, University of Houston, Houston, TX, 77204-4005, USA
[2] Midcontinent Independent System Operator (MISO), P.O. Box 4202, Carmel, IN 46082-4202, USA
[3] Eaton Corporation, Pune, Maharashtra, 411014, India
[*] xli82@uh.edu

**Abstract:** Transmission networks consist of thousands of branches for large-scale real power systems. They are built with a high degree of redundancy for reliability concern. Thus, it is very likely that there exist various network topologies that can deliver continuous power supply to consumers. The optimal transmission network topology could be very different for different system conditions. Transmission network topology control (TNTC) can provide the operator with an additional option to manage network congestion, reduce losses, relieve violation, and achieve cost saving. This paper examines the benefits of TNTC in reducing post-contingency overloads that are identified by real-time contingency analysis (RTCA). The procedure of RTCA with TNTC is presented and two algorithms are proposed to determine the candidate switching solutions. Both algorithms use available system data: sensitivity factors or shifting factors. The proposed two TNTC approaches are based on the transmission switching distribution factor (TSDF) and flow transfer distribution factor (FTDF) respectively. FTDF based TNTC approach is an enhanced version of TSDF based TNTC approach by considering network flow distribution. Numerical simulations demonstrate that both methods can effectively relieve flow violations and FTDF outperforms TSDF.

## 1. Introduction

A power system is an electrical network consisting of various elements that are used to generate, transmit, and consume electric power. Every power system contains four major components: power generation, meshed transmission network, radial distribution network, and power consumption.

It is neither cost-effective nor practical to store electric energy on a large scale. Therefore, electricity supply must meet electricity demand simultaneously for power systems. This creates serious challenges for power system real-time operations. Uncertainties associated with volatile demand and intermittent renewable generation, along with system physical constraints and reliability requirements, make real-time operations even more complex. Thus, system operators use computer-aided tools including supervisory control and data acquisition (SCADA) and energy management system (EMS) to monitor, control, and optimize the system. The measurement data from millions of measuring devices will be collected and processed by SCADA. Then, those data will be forwarded to EMS for in-depth analysis. Real-time contingency analysis (RTCA) is a key module of EMS; it scans the power system and identifies potential vulnerabilities. The RTCA results will be used by a subsequent EMS module, real-time security-constrained economic dispatch (RT SCED), to create generation redispatch solutions to eliminate the vulnerabilities.

The topology of the transmission network is traditionally considered to be fixed in EMS unless unexpected outage events occur. However, prior efforts demonstrate utilizing flexible transmission network topology control (TNTC) can achieve various significant benefits. It is demonstrated in [1] that including TNTC into optimization based scheduling and dispatching problems will increase the feasible set, which leads to a solution that is at least as good as the solution for the same problem without TNTC. In addition, operators are allowed to change transmission network topology in practice; multiple independent system operators (ISOs) acknowledge the effectiveness of TNTC and consider TNTC as a control mechanism [2]-[6]. For example, a California ISO's report states that a switching action relieved the network congestion caused by outages in the transmission system [2]; TNTC is listed as one of the effective corrective actions in the PJM's manual to relieve thermal overloaded transmission [3]-[4].

Renewable generation such as wind generation has been growing significantly in recent years. This will reduce carbon emission, relieve global warming, and contribute to a sustainable electric energy system. However, this also introduces additional uncertainty with a substantial degree into the power grid, which imposes new challenges to system reliability and incurs extra system operation cost. One technology to facilitate the grid integration of renewable generation is transmission network topology control. It shows the proposed chance-constrained transmission network topology optimization can accommodate higher utilization of wind power [7]. The chance constraints used in this optimization model can guarantee a minimal usage of wind energy at a certain probability. This work also shows TNTC can reduce generation cost of thermal generators for a system with high penetration of wind power. In [8], TNTC is used as a recourse action in the day-ahead unit commitment problem when large-scale renewable generation exists in the power system. It shows the cost can be reduced by 3.3% with the proposed TNTC resource for heavy load and large-scale generation situations. An alternative option to manage intermittent renewable generation is to set do-not-exceed (DNE) limits for intermittent renewable generations. The work in [9] illustrates that TNTC can increase the DNE limits by over 20% and lower the operation cost by more than 6%.

A zonal approach that can efficiently determine DNE limits is examined in [10].

Network congestion is a key factor that affects the system operation efficiency. The annual total congestion cost in the PJM system is over 600 million dollars in 2013 and it increases to over 1,000 million dollars in 2017 [11]. Since TNTC is a power flow control technology, it can be used to manage congestion and improve efficiency. Reference [12] studies the problem of optimal network configuration for congestion management and shows that TNTC can alleviate overloads and avoid expensive generation or load curtailment. Two deterministic and genetic based approaches are proposed and tested against several cases. Both methods show TNTC can completely relieve all overloads. It is shown in [13]-[14] that TNTC can also effectively relieve post-contingency congestion and substantially reduce congestion cost as a corrective mechanism. Although TNTC may be subject to several practical factors, it is a realistic mechanism; for instance, [15] shows that TNTC is able to relieve congestion even when the generator rotor shaft impact caused by switching is considered.

Due to constant demand increase and dramatic change in generating resources, power system expansion planning is key to ensuring future electric energy supply to consumers. Considering TNTC in power system expansion planning can defer the investment of new elements such as transmission or generator, lower the capacity of new resources, and remove the need of equipment upgrade, which would substantially reduce the cost and improve power system operation efficiency [16]. A two-stage stochastic programming model is proposed to co-optimize TNTC and investment of transmission and energy storage system [17]. Simulation results in [17] show the incorporation of TNTC in the stochastic model can reduce the total cost, including new line investment cost, new storage unit investment cost and expected operational cost, by 17% and reduce the total capacity of energy storage system by 50%. Including TNTC in system expansion planning could make the mixed integer programming problem much more complex; however, Benders decomposition can accelerate the solving process and provide quality solutions [18].

Optimization model based power system scheduling and dispatching applications provide least cost operational solutions for generating resources that meet reliability requirements to ensure continuous power supply to consumers. Including TNTC in those models may significantly enhance system security and reduce the total cost that may include the startup cost, no-load cost and operational cost. It is demonstrated that TNTC can achieve great enhancements in service reliability [19]. The cost saving with TNTC is shown to be 25% and 15% without and with *N*-1 reliability constraints on the IEEE 118-bus system [20]-[21]. With robust optimization, TNTC solutions would be feasible for a range of system operating states [22]. Short-circuit current is one major concern in power system operations; the studies conducted in [23] provide new insights to reduce the short-circuit current with TNTC. As a corrective action, TNTC can substantially reduce post-contingency violations and the proposed vicinity-base local search and data mining heuristics can provide multiple quality solutions [24]-[29]. Authentic ISO New England data and software are used to demonstrate the benefits of using TNTC as a corrective mechanism in response to both the *N*-1 and *N*-1-1 events: corrective TNTC can alleviate thermal overloads and achieve economic benefits [30]. Numerical simulations with actual market data and in-house market software show that TNTC can improve system reliability and save millions of dollars each year for ISO New England. Moreover, TNTC can address multi-element contingency. The proposed TNTC algorithm in [31] can provide effective solutions that improve the deliverability of reserves and quickly restore the loads. Furthermore, it is shown that undesired renewable generation curtailment can be avoided with TNTC by relieving post-contingency network congestion [32].

Improving power system resilience is key to withstanding high impact low probability events. Properly addressing system vulnerabilities under extreme events can effectively prevent wide area and long duration outages, enhance public safety and health, and ensure social welfare. TNTC is one of the techniques to improve power system resilience. The quantification of power system resiliency is studied in [33], which suggests TNTC is an effective approach for resiliency improvement. TNTC is considered as a temporary operation mechanism for outage recovery with minimum additional cost. TNTC can also be incorporated in the splitting strategy to enhance power grid resilience [34]; the proposed splitting strategy is used as the last resort before a critical transition occurs.

In this work, we apply TNTC in the application of RTCA. RTCA identifies potential system vulnerabilities in the form of post-contingency violations. These violations can be relieved or even eliminated by flexible topology control. In the post-contingency situations with one or multiple branches overloaded, switching a branch offline will change the network topology and then change the flow distribution in the transmission network. Beneficial switching lines may relieve the violations by reducing the flows on overloaded lines [25]. Since the contingency-induced overloads can be eliminated or reduced by corrective TNTC, the monitored branch set in RT SCED could be reduced, and higher values could be used as the limits for the associated monitored branches. Thus, the network constraints that are included in RT SCED could be relaxed; this would reduce the need for expensive generation redispatch and thus save cost [13]-[14]. Specifically, this work will examine how much violation reduction can be achieved with corrective TNTC that redistributes the network flows via switching actions. The beneficial switching actions can be determined by the proposed sensitivity factor based TNTC approaches. The proposed sensitivity factors include: transmission switching distribution factor (TSDF) and flow transfer distribution factor (FTDF). Case studies demonstrate the effectiveness of these two TSDF and FTDF based TNTC approaches.

The contributions of this paper are summarized as follows:
- The proposed TSDF and FTDF can be calculated very fast as they are derived from the simplified linearized DC power flow model; in addition, practically available sensitivity factors can be leveraged to determine TSDF and FTDF.
- The proposed TSDF and FTDF based TNTC approaches are very effective as they can efficiently identify the beneficial TNTC switching solutions to relieve post-contingency overloads in an AC setting.

- Numerical simulations on the IEEE 24-bus system in an AC setting demonstrate the effectiveness of the proposed TNTC approaches; and numerical simulations on the large-scale 2383-bus Polish system in an AC setting show the scalability of the proposed TNTC algorithms.

The rest of this paper is organized as follows: Section II presents the methodology for using topology control in RTCA. Section III describes the proposed sensitivity factors and TNTC algorithms. Case studies are presented in Section IV. Finally, Section V concludes the paper.

## 2. Methodology for RTCA with TNTC

This section presents the proposed sensitivity factor based methodology for transmission network topology control that aims to address flow violation caused by generator contingency or line contingency. Figure 1 shows the flowchart of RTCA with TNTC and the detailed steps are presented below.

- Step 1: Monitor system status and perform base-case AC power flow.
- Step 2: Determine the contingency list.
- Step 3: Conduct AC real-time contingency analysis.
- Step 4: Select critical contingencies.
- Step 5: Determine the list of candidate switching solutions for a critical contingency $c$.
- Step 6: Evaluate the TNTC candidate switching actions and identify beneficial solutions in an AC setting.
- Step 7: Set $c=c+1$ and continue the TNTC process until all critical contingencies are examined.

Power system real-time operations mainly comprise system monitoring, real-time contingency analysis, and security-constrained economic dispatch. Power system control centre collects measurements from remote terminal units and/or local control centres. After all measurements are collected, state estimation will be performed to estimate the system status and provide a basis for all other subsequent applications. Then, a base-case AC power flow simulation will be conducted to determine the solution that can accurately represent the system status (Step 1). Though state estimation and power flow are both to obtain the system base case, state estimation can be considered as a relaxed version of the power flow problem; in other words, the power flow solution is more reliable and robust than state estimation solution as state estimation solution is used as a starting point for power flow algorithms. Thus, the power flow solution will be used as the starting point for AC contingency simulations.

Step 2 is to determine a list of candidate contingencies. Due to time restriction for RTCA, a real-time application, it might not be very practical to simulate all possible contingencies. However, to validate the performance of the proposed sensitivity factor based transmission network topology control approach, the contingency list used in this work compromises all possible generator contingencies and line contingencies excluding radial line contingencies. A radial line contingency that leads to network isolation is excluded in the candidate contingency list, which is consistent with industrial practices. Radial line contingencies are addressed with other special protection schemes, which is beyond the scope of this work.

After the contingency list is determined, system operators will conduct RTCA simulating all the contingencies in that list and record the simulation results (Step 3). The recorded information includes contingency element, violation type and element, absolute violation magnitude, and relative violation magnitude in percent. Since this work focuses on line overload reduction, the contingencies that cause flow violations beyond the line short-term emergency rating will be selected as critical contingencies (Step 4).

Step 5 will first calculate the required sensitivity factors for a critical contingency identified in Step 4. Based on the sensitivity factors, TNTC algorithms will conduct and determine the candidate switching list. Details about the proposed sensitivity factors and candidate line selection method will be presented in Section III in details. Then, for the given critical contingency, Step 6 will enumerate all the switching actions in the candidate list to identify the top five beneficial TNTC solutions that determine effective network topologies for reducing contingency-induced flow violations. This process iterates until all critical contingencies are examined (Step 7).

This work can accurately reflect the flow change due to contingency events and TNTC switching actions since full AC power flow model is used. In an AC setting, the two types of state variables are the voltage magnitude and phase angle that form the basis for nodal power injection equations and branch flow equations. Fig. 2 shows a typical branch model. The branch current $I_{ij}$ and complex power $S_{ij}$ can be expressed in (1) and (2).

$$I_{ij} = y_{ij}(V_i - V_j) + j\frac{b_k}{2}V_i \quad (1)$$

$$S_{ij} = V_i I_{ij}^* \quad (2)$$

where $y_{ij}$ denotes the branch admittance; $b_k$ denotes the branch charging shunt susceptance; $V_i$ denotes the voltage at bus $i$; and $I_{ij}^*$ denotes the conjugate of branch current.

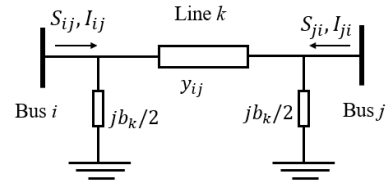

**Fig. 2.** Diagram of a typical branch model

By substituting (1) into (2), we can obtain the expressions for both active power flow $P_{ij}$ and reactive power flow $Q_{ij}$:

$$P_{ij} = V_i^2 g_{ij} - V_i V_j(g_{ij}\cos\theta_{ij} + b_{ij}\sin\theta_{ij}) \quad (3)$$

$$Q_{ij} = -V_i^2\left(b_{ij} + \frac{b_k}{2}\right) + V_i V_j(b_{ij}\cos\theta_{ij} - g_{ij}\sin\theta_{ij}) \quad (4)$$

where $g_{ij}$ and $b_{ij}$ denote the real part and imaginary part of the branch admittance $y_{ij}$.

An accurate power flow solution meets nodal power balance equations (5) and (6) within a prespecified tolerance. In other words, the values of the left-hand sides in (5) and (6) with the solved voltage magnitude and phase angle should be less than a small number that represents the desired accuracy.

$$P_i^{Inj} - V_i \sum_{j\in\{i,1,...,N\}} V_j(G_{ij}\cos\theta_{ij} + B_{ij}\sin\theta_{ij}) = 0 \quad (5)$$

$$Q_i^{Inj} - V_i \sum_{j\in\{i,1,...,N\}} V_j(-B_{ij}\cos\theta_{ij} + G_{ij}\sin\theta_{ij}) = 0 \quad (6)$$

where $G_{ij}$ and $B_{ij}$ represent the conductance and susceptance elements in the bus admittance matrix; $\theta_{ij}$ denotes the phase angle difference between bus $i$ and bus $j$; $P_i^{Inj}$ and $Q_i^{Inj}$ denote the active and reactive power injected from bus $i$.

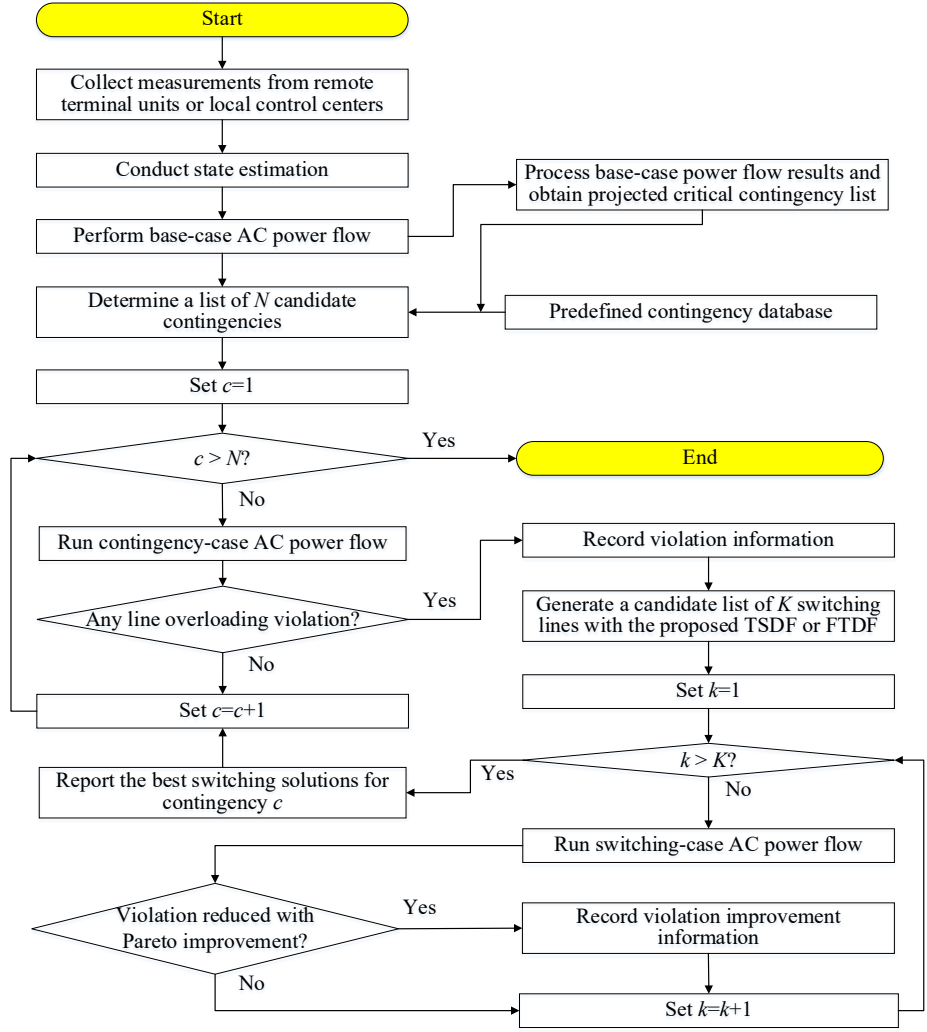

*Fig. 1.* The flowchart of TNTC-based RTCA

A stable system status should be within system physical limits, including the bus voltage limits as shown in (7) and the branch flow limits as shown in (8); otherwise, undesired system violations would be observed. This work focuses on branch flow limit check.

$$V_{i,min} \leq V_i \leq V_{i,max} \quad (7)$$

$$-Rating_i \leq \{S_{ij}, S_{ji}\} \leq Rating_i \quad (8)$$

### 3. Sensitivity Factor based TNTC Algorithms

How to efficiently determine the candidate switching list is a key question for utilizing TNTC to relieve flow violations. Most prior methods are based on optimization techniques, which are very slow and may not meet the time requirement for RTCA, a real-time application. This paper addresses this concern by proposing two fast TNTC algorithms that use the transmission switching distribution factor and flow transfer distribution factor to identify candidate switching list. Neither TSDF nor FTDF based TNTC method will include radial lines in the candidate switching list because disconnection of radial lines will cause network separation and may result in serious transient stability issue and multi-area power imbalance issue.

Several well-known and widely used sensitivity factors include power transfer distribution factors (PTDF) [35] and line outage distribution factors (LODF) [36]. PTDFs are used to model branch flow change due to nodal power injection change while LODFs model branch flow change due to line outage. The following two equations explain the meaning of PTDF and LODF respectively.

$$PTDF_{l,i} = P_l - P_l^o \quad (9)$$

$$LODF_{l,c} = \frac{P_{lc}-P_{l0}}{P_{c0}} \quad (10)$$

where $PTDF_{l,i}$ denotes the incremental flow change in real power on line $l$ due to 1 MW injection at bus $i$ and 1 MW withdrawal at the slack bus; $P_l$ denotes the flow on line $l$ after this change while $P_l^o$ denotes the original flow before this change. In (10), $LODF_{l,c}$ denotes the ratio of change in real power flow on line $l$ (due to line $c$ outage) to the pre-outage real power flow on line $c$; $P_{lc}$ denotes the flow on the monitored line $l$ after the outage of line $c$ while $P_{l0}$ and $P_{c0}$ denotes the original flows on line $l$ and line $c$ before the outage of line $c$ respectively.

The proposed TSDFs represent the ratio of change in power flow on a monitor line (due to a line switching action under contingency $c$) to the pre-switching flow on the identified switching line. TSDFs involves three elements: the contingency element, the switching line, and the post-contingency overloaded line. Specifically, as shown in (11), $TSDF_{m,k}^c$ indicates flow changes on overloaded line $m$ due to switching action on line $k$ in the system without the contingency element $c$. TSDF can be calculated with (12) by using the PTDFs for the post-contingency network without the contingency element $c$. TSDF is very similar to post-contingency LODF but with a focus on the impact of line switching on overloaded branches, which can substantially save computing time.

Unlike PTDF, LODF and TSDF that are purely dependent on the network topology and line impedance, the proposed FTDF also captures the impact of the magnitude of switching line flow. $FTDF_{m,k}^c$ can be calculated with (13) by multiplying $TSDF_{m,k}^c$ with switching line flow $P_{k,c}$ under contingency $c$. For PTDF, LODF and TSDF, the elements are within the range of [-1, 1] unless there exist lines with negative reactance; for instance, line over compensation may lead to equivalent negative reactance, which would result in PTDF values greater than 1 or smaller than -1. As for FTDF, there is no such limit since the flow magnitude on the switching line under contingency may vary substantially.

$$TSDF_{m,k}^c = \frac{\Delta P_{m,c,k}}{P_{k,c}} = \frac{P_{m,c,k} - P_{m,c}}{P_{k,c}} \quad (11)$$

$$TSDF_{m,k}^c = \frac{PTDF_{m,f(k)}^c - PTDF_{m,t(k)}^c}{1 - (PTDF_{k,f(k)}^c - PTDF_{k,t(k)}^c)} \quad (12)$$

$$FTDF_{m,k}^c = TSDF_{m,k}^c P_{k,c} \quad (13)$$

where $c$ denotes the contingency line; $m$ denotes the overloaded line due to contingency $c$; $k$ represents the switching line; $P_{k,c}$ denotes the flow on switching line $k$ under contingency $c$; $P_{m,c}$ denotes the flow on the overloaded line $m$ under contingency $c$; $P_{m,c,k}$ denotes the flow on line $m$ after line $k$ is switched out of service under contingency $c$; $\Delta P_{m,c,k}$ denotes the flow change on line $m$ due to line $k$ switched out of service under contingency $c$; $f(k)$ denotes the from-bus of line $k$ while $t(k)$ denotes the to-bus of line $k$; $PTDF_{m,n}^c$ denotes the PTDF for flow change on line $m$ due to injection change at node $n$ for the network without contingency element $c$.

Since the proposed TNTC scheme is to handle flow violations caused by contingencies in real-time, computing time is a key performance indicator. TNTC switching simulations are essentially a series of AC power flow simulations with different network topologies. As a result, the algorithm computing time depends on the number of candidate switching lines. Thus, it is very important to select as least candidate switching lines as possible to meet time requirements. On the other hand, it may not find any beneficial switching solutions with a very small set of candidate TNTC switching lines. Thus, to determine a proper number of TNTC switching candidates, simulations with different numbers of candidate TNTC switching lines are conducted; they are referred to as FTDF'$N$' or TSDF'$N$' which indicates a list of $N$ candidate switching lines.

Once the sensitivity factors, either TSDF or FTDF, are available, a ranking method can be used for candidate line selection. If the actual flow direction of the overloaded line is the same with its reference direction, the candidate switching lines for relieving the overload on this line will be ordered from small values of TSDF or FTDF to large values of TSDF or FTDF. For instance, FTDF20 will select 20 lines that correspond to the 20 smallest FTDF factors (negative numbers) for the given overloaded lines. However, if the actual flow direction of the overloaded line is different than its reference direction, then the candidate switching list determined by FTDF20 would consist of lines that have the 20 largest values of FTDF factors (positive numbers) for the given overloaded lines.

In this paper, we use TSDF and FTDF to determine the candidate switching lines. However, these factors are based on the simplified linearized DC power flow model. It cannot guarantee that the preselected candidate switching lines will perform as expected; even AC feasibility cannot be guaranteed. Thus, we need to evaluate those candidate switching lines with the full AC power flow model. This is indicated in Figure 1.

Enumerating all switchable lines for each critical contingency will ensure the best TNTC switching solution but it comes at the cost of a much longer computing time. However, complete enumeration (CE) is implemented in this work to evaluate the performance of the proposed TNTC algorithms.

Pareto improvement is enforced throughout this paper. A TNTC solution with Pareto improvement will reduce the total violation while not increasing any existing individual violations or creating new violations. Pareto improvement is important because no individual entity accepts violations to their own facilities to relieve overloading issue on the lines belonging to another entity. To measure the performance of a single TNTC switching solution, violation reduction in percent (VRP) for switching line $k$ against contingency $c$ is proposed and defined in (14). To measure the performance of TNTC algorithms, the average percent violation reduction $\varepsilon$ is defined in (15) as the average of violation reduction in percent among all contingencies. Depth is used to evaluate the effectiveness of the proposed TNTC algorithms. It is defined as the rank of the identified TNTC solution in the candidate switching list.

$$VRP_{m,k}^c = \frac{\gamma_m^c - \gamma_{m,k}^c}{\gamma_m^c} \quad (14)$$

$$\varepsilon = \frac{1}{N_C} \sum_{c \in C} \frac{\gamma_m^c - \gamma_{m,k^*}^c}{\gamma_m^c} \quad (15)$$

where $\gamma_m^c$ denotes the violation on line $m$ under contingency $c$; $\gamma_{m,k}^c$ denotes the violation on line $m$ after switching line $k$ out of service under contingency $c$; $\gamma_{m,k^*}^c$ denotes the violation on line $m$ with the best switching action $k^*$ for contingency $c$; set $C$ denotes the set of critical contingencies and $N_C$ denotes the number of critical contingencies.

The concept of TNTC for relieving branch overloads is counterintuitive since TNTC switches off a second branch after the system already losses a branch. The disconnection of a branch might increase the overall equivalent impedance from one area to another area and make the transmission network weaker in some situations; however, it is not always the case, and it does not mean the transfer capacity between areas will decrease accordingly. The transmission network is highly meshed and different branches have different impedances and capacities; the power flow distribution in the transmission network must follow physical rules. It is very common that some branches are overloaded while their parallel branches still have a lot of available capacity left. The switching of a specific branch may transfer some flow on the overloaded lines to the parallel lines. This will reduce the overloads; and the loading level of non-overloaded parallel lines may increase but it may not exceed or reach the limits. Thus, TNTC may relieve network congestion and improve power transfer capacity between different areas. In addition, the switching line is not a random line; it is identified by the proposed verifiable approach that includes a step to validate the TNTC solution by evaluating the post-switching system condition via switching-case AC power flow simulations. The positive effect of TNTC would be more obvious if the grid is under a more stressed situation, since more congestions exist in heavily loaded system and there would be a higher chance of finding beneficial switching solutions. Thus, the proposed TNTC approach can be considered for adoption by any practical power systems.

It is worth noting that there is a possibility that switching on lines that were previously switched off may also relieve overloads. Since the candidate lines for switching on are much less than the candidate lines for switching off in practice, this paper focuses on the study of switching off branches. TNTC also has the potential to mitigate voltage violation; for instance, switching off a heavily loaded transmission line that consumes a lot of reactive power may help improve under-voltage condition. Since TSDF and FTDF are derived from the linearized DC power flow model, they provide no information regarding voltage and reactive power; thus, this paper focuses on branch flow violation relief only.

## 4. Numerical Simulations

In this section, two cases are tested to demonstrate the effectiveness of TNTC in terms of violation relief. The small-scale test case is the IEEE 24-bus test system - one area of the RTS96 system [37]; this case is mainly used to illustrate the fundamental of TNTC and how it reroutes power in the network to achieve violation relief. The large-scale test case is the Polish system with over 2,000 buses, which is used to demonstrate the proposed TNTC algorithms and showcase the algorithm scalability.

### A. IEEE 24-bus test system

The IEEE 24-bus system has 24 buses and 38 branches. The network topology of this system with bus number and branch number is shown in Fig. 3 [38]. It has two zones: the 138-kV zone consisting of bus 1 through bus 10 and the 230-kV zone consisting of bus 11 through bus 24. These two zones are connected by five transformers. In the 138-kV zone, the total load is 1,137 MW while the total generation is 512 MW. This indicates that the 138-kV zone is lack of generation and imports power from the 230-kV zone through those five transformers. The original IEEE 24-bus system is a highly redundant system with high thermal limits and no congestions would be observed. Thus, in order to show the effect of TNTC, branch thermal limits are reduced in the revised IEEE 24-bus system used in this paper. For example, the long-term thermal limits for branch 23, 24, 25, 26 are all set to 240 MVA while their short-term emergency thermal limits are all set to 275 MVA in the revised IEEE 24-bus system.

RTCA is conducted on 37 non-radial branches and identifies two critical contingencies: (i) branch 7 connecting bus 3 and bus 24; (ii) branch 27 connecting bus 15 and bus 24. Since branch 7 and branch 27 are on the same path connecting the 230-kV zone and the 138-kV zone, the impact of branch 7 outage and branch 27 outage would be very similar. The outage of branch 7 results in a single violation of 26.4 MVA on branch 23 while the outage of branch 27 leads to a single violation of 26.3 MVA on branch 23. The failure of a tie-path will transfer the power it carries to other tie-paths, which may overload one or multiple branches on the remaining tie-paths.

After two critical contingencies are identified, three methods including CE, TSDF20 and FTDF20 are conducted. CE considers all non-radial lines as candidate switching solutions while the candidate switching lists for both TSDF20 and FTDF20 consist of 20 non-radial lines that are determined by the proposed TSDF and FTDF based methods respectively. The violation reduction in percent with these methods under the contingency on branch 7 is illustrated in Fig. 4. The x-axis shows the order of top five best TNTC solutions. It is observed that FTDF20 can achieve very similar results with CE and outperforms TSDF20. The TNTC time for determining the optimal switching solutions to the two critical contingencies with FTDF20 is only 0.031 seconds. More detailed results for FTDF20 are shown in Table 1. The selected TNTC solutions are to balance the flow distribution on the remaining tie-paths so that no single tie-path is overloaded. For instance, by switching branch 16 out of service in the post-contingency situation, less power will go through the path containing branch 23 and branch 19, which helps substantially reduce the post-contingency flow violation on branch 23.

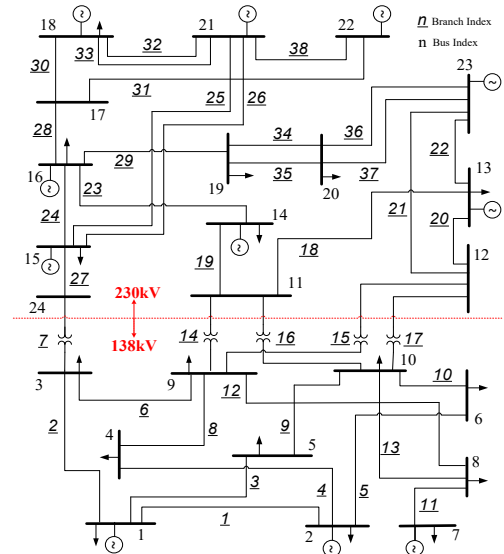

*Fig. 3.* The network topology of the IEEE 24-bus system [38]

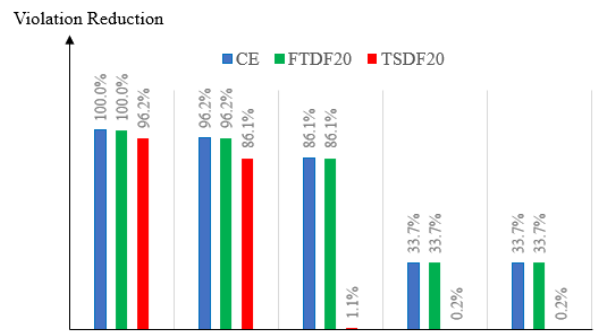

*Fig. 4.* Violation reduction with different TNTC methods

**Table 1** TNTC Results with Pareto improvement on IEEE 24-bus system

| Contingency branch | Violation reduction in percent with TNTC solutions | | | | |
|---|---|---|---|---|---|
| | 1st Best | 2nd Best | 3rd Best | 4th Best | 5th Best |
| 7 | 100% (Brc:19) | 96.2% (Brc:16) | 86.1% (Brc:14) | 33.7% (Brc:36) | 33.7% (Brc:37) |
| 27 | 100% (Brc:19) | 96.4% (Brc:16) | 86.4% (Brc:14) | 33.8% (Brc:36) | 33.8% (Brc:37) |

"Brc" denotes branch; the information in the parenthesis denotes the associated TNTC switching solution.

A complete example to show the effectiveness of the proposed TNTC approach is presented as follows:
- In the base case (pre-contingency situation), the flow on branch 23 is 213.6 MVA, which is below the long-term normal thermal limit of 240 MVA. There are no violations through the entire system.
- After the outage of branch 27, the flow on branch 23 becomes 301.3 MVA, which is well above the short-term emergency thermal limit of 275 MVA. In other words, the outage of branch 27 results in a violation of 26.3 MVA on branch 23. There are no other violations in the network.
- After an overload of 26.3 MVA is observed in the post-contingency situation, TNTC is conducted and it identifies a beneficial switching solution: line 19. After switching line 19 off, the flow on branch 23 reduces to 138 MVA that is well below its thermal limit and no other lines are overloaded. Thus, the identified TNTC solution can fully eliminate the overload.

### B. The Polish System

To validate the proposed TNTC approach and show the scalability of the proposed switching algorithms, the large-scale Polish system is tested. This system has 2,383 buses that are connected by 2,895 branches including 2,724 transmission lines and 171 two-winding transformers. The total load for active power and reactive power is 21.4 GW and 6.6 GVAr respectively.

RTCA is conducted on this Polish case to identify the critical contingencies that would lead to branch violations. It takes 11.6 seconds to simulate 327 generator contingencies and 79.8 seconds to simulate 2251 non-radial branch contingencies. Three generator contingencies and seven branch contingencies are identified to be critical contingencies.

Table 2 shows the statistics of the critical contingencies identified by RTCA. The average violation per contingency is 11.8 MVA with a standard deviation of 7.9 MVA. The maximum violation on a single line is 30 MVA, which corresponds to 10.6% of the associated line emergency rating. The minimum violation is 5 MVA, which corresponds to 4.9% of the associated line emergency rating. The detailed RTCA results are shown in Table 3. As shown in Table 3, most critical contingencies result in a single overload only while one branch contingency causes two overloads.

**Table 2** Statistics about the critical contingencies identified by RTCA

| Number of critical contingencies | Violation (MVA) | | | | |
|---|---|---|---|---|---|
| | Maximum | Minimum | Average | Median | Standard deviation |
| 10 | 30.0 | 5.0 | 11.8 | 8.5 | 7.9 |

Table 4 presents the violation relief results with the first best TNTC solutions for various algorithms. For branch contingency 7, TSDF15 and TSDF20 algorithms can identify the best TNTC solution that reduce the violation by 30.2% since its ranking is 14 in the TSDF-based candidate list; the best TNTC solutions for FTDF15 and FTDF20 can only reduce the violation by 4.4% and 20.1% respectively. However, for the rest nine critical contingencies, the FTDF algorithm with a list of only five candidates can achieve the same results with CE. TSDF10 can achieve the same results with CE for six critical contingencies and TSDF20 can achieve the same results with CE for eight critical contingencies. TSDF methods struggle to locate the best TNTC solutions for generator contingency 1, generator contingency 2, branch contingency 3, and branch contingency 7 while FTDF methods only struggle to identify the best solution for branch contingency 7. The main reason for why FTDF outperforms TSDF is that FTDF considers the amount of flow on switching lines while TSDF does not.

**Table 3** RTCA Results on the Polish system

| Contingency type | Generator Contingency | | | Branch Contingency | | | | | | |
|---|---|---|---|---|---|---|---|---|---|---|
| Contingency index | 1 | 2 | 3 | 1 | 2 | 3 | 4 | 5 | 6 | 7 |
| Flow violation (MVA) | 30.0 | 14.5 | 9.6 | 5.0 | 7.7 | 21.0 | 7.4 | 8.5 | 8.5 | 6.0 |
| # of overloaded lines | 1 | 1 | 1 | 1 | 1 | 2 | 1 | 1 | 1 | 1 |

**Table 4** TNTC Results with the first best switching solutions for various algorithms

| Contingency type | Generator Contingency | | | Branch Contingency | | | | | | |
|---|---|---|---|---|---|---|---|---|---|---|
| Contingency index | 1 | 2 | 3 | 1 | 2 | 3 | 4 | 5 | 6 | 7 |
| Algorithm | Violation reduction in percent | | | | | | | | | |
| CE | 100% | 100% | 100% | 100% | 100% | 82.3% | 100% | 100% | 100% | 30.2% |
| TSDF5 | 37.6% | 75.5% | 0% | 100% | 0% | 100% | 100% | 100% | 0% | |
| TSDF10 | 41.2% | 75.5% | 100% | 100% | 5.5% | 100% | 100% | 100% | 0% | |
| TSDF15 | 41.2% | 75.5% | 100% | 100% | 5.5% | 100% | 100% | 100% | 30.2% | |
| TSDF20 | 41.2% | 75.5% | 100% | 100% | 82.3% | 100% | 100% | 100% | 30.2% | |
| FTDF5 | 100% | 100% | 100% | 100% | 100% | 82.3% | 100% | 100% | 100% | 0% |
| FTDF10 | 100% | 100% | 100% | 100% | 100% | 82.3% | 100% | 100% | 100% | 0% |
| FTDF15 | 100% | 100% | 100% | 100% | 100% | 82.3% | 100% | 100% | 100% | 4.4% |
| FTDF20 | 100% | 100% | 100% | 100% | 100% | 82.3% | 100% | 100% | 100% | 20.1% |

Table 5 shows the statistics of the TNTC results. The average violation reduction in percent ($\varepsilon$) for FTDF10 and FTDF20 are 88.2% and 90.2% respectively; as for TSDF10 and TSDF20, the average violation reduction in percent are only 72.2% and 82.9% respectively. The average number ($\mu$) of TNTC solutions that can fully eliminate the associated contingency-induced violations is 5.5 for FTDF10 and 6.25 for FTDF20; however, it drops to 3.7 for both TSDF10 and TSDF20. By comparing (i) $n_1$ - the number of contingencies that the associated violations are fully eliminated with TNTC, (ii) $n_2$ - the number of contingencies that the associated violations are partially reduced with TNTC, (iii) $n_3$ - the number of contingencies that TNTC does not help at all, it is obvious that the FTDF approach outperforms the TSDF approach. Fig. 5 shows the computational time for different TNTC methods. It illustrates the computational complexity for FTDF methods and TSDF methods are very similar and they are much faster than complete enumeration. The solution time for the TSDF and FTDF methods on the Polish system is below 10 seconds, which demonstrates the scalability of the proposed algorithms.

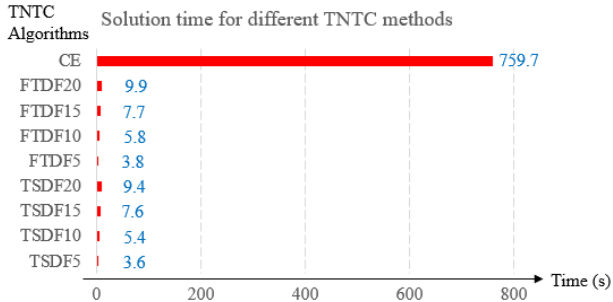

*Fig. 5.* Computational time for various TNTC methods

Table 6 shows the average violation reduction in percent over all critical contingencies with the top five best switching solutions. The results obtained with the algorithms FTDF10, FTDF15, and FTDF20 are very similar to the results of CE for any of the top five best TNTC solutions. Even with the fifth best switching solutions of FTDF10, the average violation reduction in percent is 79.9%, which indicates that multiple beneficial TNTC solutions exist and they can substantially reduce the contingency-induced violations with Pareto improvement. Though the first best TNTC solutions provided by FTDF5 can achieve similar results with FTDF10, FTDF15, FTDF20, and CE, the performance of the second, third, and fourth best TNTC solutions of FTDF5 are 10%, 18%, and 26% lower than CE respectively. Even worse, the fifth best TNTC solution of FTDF5 can only provide a violation reduction of 15% that is much lower than other FTDF algorithms. The main reason of significantly reduced performance of FTDF5 is that its candidate switching list is too short and some important beneficial TNTC solutions are not included in this candidate list. To conclude, a candidate list comprising 10 switching actions is required to find five best switching solutions. As for TSDF5, TSDF10, TSDF15 and TSDF20, their performance is obviously lower than the FTDF methods. Even with 20 lines in the candidate list for TSDF20, the violation reduction is much lower than FTDF10 that has a candidate list of 10 lines only. This is because unlike FTDF, TSDF methods do not consider the flow on the switching line in the post-contingency situation that may significantly affect the results.

The total violation caused by the 10 critical contingencies are 118.2 MVA, which can be substantially reduced with TNTC. The total violations in the post-switching situations with the five best TNTC switching solutions are listed in Table 7. With the top TNTC solutions reported by CE, the total violation is reduced to 7.9 MVA only. Method FTDF10 reduces the total violation to 9.7 MVA, which is only 2 MVA higher than CE. TSDF20 corresponds to a total violation of 29 MVA, which is three times higher than CE and FTDF10; TSDF10 has a total violation of 47 MVA which is about five times higher than FTDF10. By comparing the total violations in the post-switching situations, it is demonstrated that though TSDF methods can achieve significant violation reduction, their performance is not as good as FTDF methods. Fig. 6 shows the system total violations with the first best TNTC solutions in the post-switching situations. Obviously, FTDF methods outperform TSDF methods; and FTDF methods can achieve very similar results with complete enumeration.

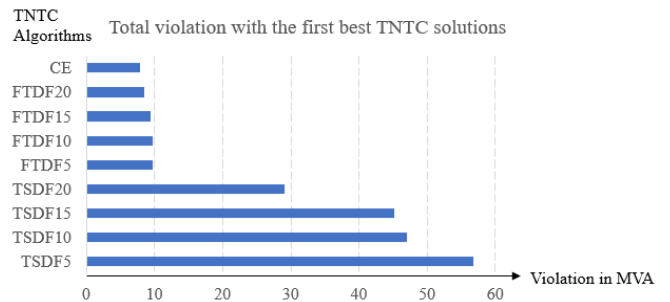

*Fig. 6.* Total overload after implementing the first best TNTC solutions

**Table 5** Statistics of TNTC results

|  | TSDF5 | TSDF10 | TSDF15 | TSDF20 | FTDF5 | FTDF10 | FTDF15 | FTDF20 | CE |
|---|---|---|---|---|---|---|---|---|---|
| Solution time (s) | 3.6 | 5.4 | 7.6 | 9.4 | 3.8 | 5.8 | 7.7 | 9.9 | 759.7 |
| $n_1$ | 5 | 6 | 6 | 6 | 8 | 8 | 8 | 8 | 8 |
| $n_2$ | 3 | 3 | 4 | 4 | 1 | 1 | 2 | 2 | 2 |
| $n_3$ | 2 | 1 | 0 | 0 | 1 | 1 | 0 | 0 | 0 |
| $\mu$ | 3.0 | 3.7 | 3.7 | 3.7 | 2.8 | 5.5 | 6.25 | 6.25 | 6.4 |
| $\varepsilon$ | 61.3% | 72.2% | 75.2% | 82.9% | 88.2% | 88.2% | 88.7% | 90.2% | 91.3% |

**Table 6** Average violation reduction in percent with the five best TNTC solutions

| Average Violation Reduction | TSDF5 | TSDF10 | TSDF15 | TSDF20 | FTDF5 | FTDF10 | FTDF15 | FTDF20 | CE |
|---|---|---|---|---|---|---|---|---|---|
| 1st Best TNTC | 61.3% | 72.2% | 75.2% | 82.9% | 88.2% | 88.2% | 88.7% | 90.2% | 91.3% |
| 2nd Best TNTC | 50.0% | 69.6% | 69.6% | 70.0% | 76.2% | 86.2% | 86.2% | 86.6% | 88.2% |
| 3rd Best TNTC | 47.5% | 64.9% | 65.9% | 66.3% | 67.1% | 84.7% | 84.7% | 85.0% | 85.2% |
| 4th Best TNTC | 36.3% | 56.0% | 56.6% | 56.6% | 55.0% | 81.4% | 81.4% | 81.4% | 81.8% |
| 5th Best TNTC | 0.0% | 20.3% | 43.0% | 49.4% | 15.2% | 79.9% | 80.5% | 80.5% | 80.8% |

**Table 7** Total violation in the post-switching situations with the five best TNTC solutions

| Total Violation (MVA) | TSDF5 | TSDF10 | TSDF15 | TSDF20 | FTDF5 | FTDF10 | FTDF15 | FTDF20 | CE |
|---|---|---|---|---|---|---|---|---|---|
| 1st Best TNTC | 56.7 | 47.0 | 45.2 | 29.1 | 9.7 | 9.7 | 9.4 | 8.5 | 7.9 |
| 2nd Best TNTC | 81.1 | 51.4 | 51.4 | 50.9 | 21.7 | 14.0 | 14.0 | 13.9 | 12.8 |
| 3rd Best TNTC | 83.0 | 60.5 | 58.9 | 58.3 | 31.1 | 17.1 | 17.1 | 16.9 | 16.8 |
| 4th Best TNTC | 92.1 | 70.0 | 69.0 | 69. | 51.9 | 26.0 | 26.0 | 26.0 | 25.8 |
| 5th Best TNTC | 118.2 | 100.9 | 77.6 | 73.9 | 105.9 | 30.1 | 28.8 | 28.8 | 28.6 |

**Table 8** Average depth of beneficial TNTC solutions

| | TSDF5 | TSDF10 | TSDF15 | TSDF20 | FTDF5 | FTDF10 | FTDF15 | FTDF20 |
|---|---|---|---|---|---|---|---|---|
| 1st Best TNTC | 1.75 | 3.2 | 4.3 | 5.2 | 2.4 | 2.4 | 3.7 | 3.9 |
| 2nd Best TNTC | 2.6 | 4.3 | 4.3 | 5.7 | 3.6 | 3.9 | 3.9 | 5.0 |
| 3rd Best TNTC | 4.0 | 5.6 | 7.1 | 7.9 | 3.3 | 5.0 | 5.0 | 6.4 |
| 4th Best TNTC | 2.8 | 5.5 | 6.1 | 7.1 | 4 | 6.4 | 6.4 | 6.4 |
| 5th Best TNTC | N/A | 8.7 | 9.3 | 12.4 | 2 | 7.2 | 7.7 | 7.7 |

Table 8 shows the average depth of the beneficial switching solutions for the proposed TNTC approaches. The average depths of the fifth best TNTC solutions for FTDF10, FTDF15, and FTDF20 are 7.2, 7.7, and 7.7 respectively. This indicates that the fifth best TNTC solutions rank around 7th and 8th in the TNTC candidate list, which explains why FTDF5 has a poor performance while FTDF10 is able to find all top five beneficial TNTC solutions.

## 5. Conclusion

Practical power systems typically consist of thousands of buses and branches that form large-scale transmission networks with redundancy. Flexible utilization of the transmission network can benefit the system in various aspects. However, the transmission networks are considered as static asset in today's operational tools including energy management system. EMS is a decision support software system for real-time operations of electric power systems. The RTCA module of EMS can determine critical contingencies and the associated violations, which is key to scanning the system and identify potential system vulnerabilities. This paper investigates the impact of TNTC on relieving the potential post-contingency flow violations identified by RTCA.

Two TNTC approaches are proposed to determine the candidate switching solutions and determine the optimal topology considering system reliability concern. They are based on transmission switching distribution factor and flow transfer distribution factor respectively. These factors can be easily calculated with existing sensitivity factors such as power transfer distribution factors; thus, the implementation of the proposed methods is straightforward and requires very short solution time. Though the proposed sensitivity factors are based on the DC power flow model, case studies with full AC model based simulations show that they can effectively identify the beneficial switching lines. Numerical simulations verify that the proposed TSDF and FTDF approaches can provide TNTC solutions that achieve substantial post-contingency violation reduction and they can scale well on large-scale power systems. Thus, the proposed switching algorithms have the potential to be implemented in practice. Simulation results also show FTDF is more efficient for identifying beneficial switching solutions than TSDF. Moreover, case studies demonstrate multiple TNTC beneficial solutions would be available for overload relief.